# Structures riemanniennes $L^p$ et K-homologie

By Michel Hilsum

### Abstract

We construct analytically the signature operator for a new family of topological manifolds. This family contains the quasi-conformal manifolds and the topological manifolds modeled on germs of homeomorphisms of $\mathbf{R}^n$ possessing a derivative which is in $L^p$, with $p > \frac{1}{2}n(n+1)$. We obtain an unbounded Fredholm module which defines a class in the K-homology of the manifold, the Chern character of which is the Hirzebruch polynomial in the Pontrjagin classes of the manifold.

This generalizes previous works of N. Teleman for Lipschitz manifolds and of A. Connes, N. Teleman and D. Sullivan for quasi-conformal manifolds of even dimension [11], [5].

### Introduction

Le but de ce travail est d'étendre la construction de l'opérateur de signature d'une variété $C^\infty$ à de nouveaux exemples de variétés topologiques.

Cette généralisation a été faite par N. Teleman [11] pour une variété linéaire par morceaux, puis pour une variété lipschitzienne. Récemment, l'existence de cet opérateur a été établie pour une variété quasi-conforme de dimension paire par A. Connes, D. Sullivan et N. Teleman [5].

Pour une variété $C^\infty$ riemannienne $(V, g)$, l'opérateur de signature est l'extension autoadjointe de $D = d + d^*$, où $d$ est la différentielle extérieure de de Rham, sur l'espace de Hilbert $\Omega^*(V, g)$ des formes différentielles de carré intégrable. Le couple $(\Omega^*(V,g), D)$ détermine alors un élément de $K_n(V)$, où $n = \dim V$. Lorsque $V$ est de dimension paire, toute l'information K-homologique de $D$ est contenu dans l'opérateur $(dd^* - d^*d)(1+D^2)^{-1}$ restreint à $\Omega^m(V, g)$, où $2m = n$ [5]. Lorsque $V$ est de dimension impaire, la classe de K-homologie de $D$ est de même exactement représentée par l'opérateur $\tau d$ agissant dans $\Omega^m(V, g)$, où $\tau$ est l'involution de Hodge, et $2m+1 = n$.

Nous construisons l'analogue de ces modules de Fredholm pour les variétés quasi-conformes et les variétés qui sont munies d'une structure dérivable pour laquelle les différentielles des changements cartes sont dans $L^p$ pour $p > \frac{1}{2}n(n+1)$.

Pour une telle variété $V$, il est possible de définir des structures riemanniennes $g$ et d'y associer les espaces de Hilbert $\Omega^k(V,g)$ de formes différentielles correspondants. Ces espaces de Hilbert ne sont pas localement isomorphes à ceux construits à partir d'une structure riemannienne $C^\infty$. Cependant, à tout recouvrement $(\mathcal{O}_i)_{i \in I}$ de $V$ par des ouverts de carte, il est possible d'associer une métrique riemannienne dont l'image dans chaque $\mathcal{O}_i$ est telle que, en notant $m$ la partie entière de $\frac{n}{2}$, pour $k = m, m+1$, on ait les inclusions:

$$L^{p_k}(\mathcal{O}_i, \Lambda^k_{\mathbf{C}}(T^*\mathcal{O}_i)) \subset \Omega^k(\mathcal{O}_i, g) \subset L^{q_k}(\mathcal{O}_i, \Lambda^k_{\mathbf{C}}(T^*\mathcal{O}_i))$$

avec $p_m^{-1} + n^{-1} > q_{m+1}^{-1}$ (cf. définition 2.1).

Pour de telles métriques, la différentielle extérieure de de Rham, au sens des distributions, définit par restriction un opérateur $d$ densément défini et fermé. Le résultat principal de cet article est alors:

THÉORÈME. *Il existe $N \geq n$ ne dépendant que de $g$ tel que:*

1) *Si $n$ est impair, l'opérateur $D = \tau d$ sur $\Omega^m(V,g)$ est densément défini, autoadjoint et $D(1+D^2)^{-1} \in \mathcal{L}^{N+}$.*

2) *Si $n$ est pair, l'opérateur $D = d + d^*$ de $\Omega^m(V,g)$ dans $\Omega^{m-1}(V,g) \oplus \Omega^{m+1}(V,g)$ est densément défini, fermé, anticommute à $\tau$ et $(1+D^2)^{-\frac{1}{2}} \in \mathcal{L}^{N+}$.*

*Dans tout les cas, l'opérateur $D$ obtenu est unique à homotopie non bornée près.*

Cet opérateur détermine une classe de la K-homologie de $V$, dont le caractère de Chern est égal, *via* la dualité de Poincaré, au polynôme d'Hirzebruch en les classes de Pontryagin rationnelles de $V$.

Le plan de ce travail est le suivant: nous démontrons d'abord ce résultat dans 1) pour le cas où $V$ est le tore de dimension $n \geq 1$ et pour une métrique riemannienne vérifiant la condition ci-dessus. Ceci permet ensuite de passer au cas général dans 2). Dans 3) et 4), nous décrivons les applications en K-théorie aux variétés topologiques quasi-conformes ou bien dérivable d'ordre strictement supérieur à $\frac{1}{2}n(n+1)$.

## 1. Métriques $L^p$ sur le tore

Dans tout ce qui suit, $n$ est un entier $\geq 1$ fixé, $m$ est la partie entière de $\frac{n}{2}$ et $U$ un sous-ensemble ouvert relativement compact de $\mathbf{R}^n$. Etant donné un fibré vectoriel complexe $E$ munie d'une structure hermitienne $h$ et $p \geq 1$, on note $L^p(U, E)$ l'espace de Banach des sections mesurables $f : U \to E$ telles que $\|f\|_p = \{\int_U \|f(x)\|^p dx\}^{\frac{1}{p}} < +\infty$ où $dx$ est la mesure de Lebesgue sur $U$.

Dans ce qui suit, une *structure riemannienne $g$ sur $U$* est une fonction, mesurable pour la classe de Lebesgue, à valeurs dans l'ensemble des structures euclidiennes sur $\mathbf{R}^n$. Une telle structure détermine un espace de Hilbert $\Omega^k(U, g)$ formé des formes différentielles complexes $\omega$, mesurables pour la classe de Lebesgue, de degré $k$ et vérifiant:

$$\|\omega\|_g^2 = \int_U \lambda^k(g)(\omega, \omega) \mu_g < +\infty$$

où $\lambda^k(g)$ (resp. $\mu_g$) désigne la forme quadratique (resp. la forme volume) canoniquement associée à $g$ sur $\Lambda_{\mathbf{C}}^k(T^*U)$.

Soit $\tau : U \to \mathcal{L}(\Lambda_{\mathbf{C}}(\mathbf{R}^{n*}))$ le champ des involutions de Hodge déterminé par $g$:

$$\tau\omega = \begin{cases} i^{p(p-1)+m} * \omega & \text{si } n \text{ est pair} \\ i^{p(p+1)+m+1} * \omega & \text{si } n \text{ est impair} \end{cases}$$

où $p$ est le degré de $\omega$ et $*$ le champ d'unitaires déterminé en $x \in U$ par $* \wedge_{i \in I} e_i(x) = \varepsilon(I) \wedge_{i \notin I} e_i(x)$ où $e_i(x)$ est une base orthonormale pour $g_x$ et $\varepsilon(I) \in \{-1, 1\}$. Le champ $\tau$ détermine un opérateur unitaire de $\Omega^k(U, g)$ sur $\Omega^{n-k}(U, g)$, ce qui montre que la forme sesquilinéaire $\alpha, \beta \to \int_U \alpha \wedge \bar{\beta}$ est bien définie sur $\Omega^k(U,g) \times \Omega^{n-k}(U, g)$ et identifie $\Omega^{n-k}(U, g)$ avec le dual conjugué de $\Omega^k(U, g)$. Si $T : \Omega^k(U, g) \to \Omega^l(U, g)$ est un opérateur linéaire densément défini, nous noterons ${}^tT$ l'opérateur transposé $\Omega^{n-l}(U, g) \to \Omega^{n-k}(U, g)$ pour cette dualité, c'est à dire l'opérateur ayant pour domaine les éléments $\omega \in \Omega^{n-k}(U, g)$ tels que l'application définie sur $\text{dom}\, T$ par $\alpha \to \int_U T\alpha \wedge \bar{\omega}$ s'étende continûment à $\Omega^k(U, g)$.

*Définition* 1.1. Soit $\mathcal{R}(U)$ *l'ensemble des structures riemanniennes $g$ sur $U$ telles qu'il existe, pour $k = m, m+1$, des réels $1 \leq q_k \leq 2 \leq p_k \leq +\infty$ et $B > 1$ satisfaisant aux relations suivantes:*

(1) $$L^{p_k}(U, \Lambda_{\mathbf{C}}^k(T^*U)) \subset \Omega^k(U, g) \subset L^{q_k}(U, \Lambda_{\mathbf{C}}^k(T^*U)),$$

(2) $$\frac{1}{p_m} + \frac{1}{n} > \frac{1}{q_{m+1}},$$

(3) $$g > B.$$

*Remarque* 1.2. 1) Lorsque $n$ est impair, on peut supposer, par dualité, que $p_m^{-1} + q_{m+1}^{-1} = 1$.

2) Lorsque $n$ est pair, par dualité, les inclusions (1) entraînent:

$$L^{p_{m-1}}(U, \Lambda_{\mathbf{C}}^{m-1}(T^*U)) \subset \Omega^{m-1}(U, g) \subset L^{q_{m-1}}(U, \Lambda_{\mathbf{C}}^{m-1}(T^*U))$$

avec $p_m^{-1} + q_{m-1}^{-1} = p_{m-1}^{-1} + q_m^{-1} = 1$ et on a encore $p_{m-1}^{-1} + n^{-1} > q_m^{-1}$. Par le théorème du graphe fermé, ces inclusions d'espace de Banach sont automatiquement continues.

Nous pouvons définir la dérivée extérieure $d_U : \Omega^k(U,g) \to \Omega^{k+1}(U,g)$ pour $k = m$ en dimension impaire et $k \in \{m-1, m\}$ en dimension paire: le domaine de $d_U$ est constitué des éléments $\omega \in \Omega^k(U,g)$ qui définissent des courants dont la dérivée extérieure distributionnelle est un élément de $\Omega^{k+1}(U,g)$. Comme $\Omega^k(U,g) \subset L^1(V, \Lambda_{\mathbf{C}}^k(T^*U))$, tous ses éléments sont des courants et $\alpha \in \mathrm{dom}\, d_U$ si et seulement s'il existe $\beta \in \Omega^{k+1}(U,g)$ tel que pour tout $\gamma \in C_c^\infty(U, \Lambda^{n-k}(T^*U))$, l'égalité suivante soit satisfaite:

$$\int_U \alpha \wedge d\gamma = (-1)^{k-1} \int_U \beta \wedge \gamma.$$

Comme $C_b^\infty(U, \Lambda^k) \subset \mathrm{dom}\, d_U$ pour $k \in \{m-1, m\}$, le domaine de $d_U$ est dense. Notons enfin que pour $f \in C_c^\infty(U)$, la multiplication extérieure par le covecteur $df$ est un opérateur continu sur $\Omega^*(U,g)$ de norme plus petite que $B^{-1}\|df\|_\infty$, ce qui exprime que le commutateur densément défini $[d_U, f]$ s'étend par continuité sur $\Omega^m(U,g)$ et:

(1.3) $$\|[d_U, f]\| \leq B^{-1}\|df\|_\infty.$$

Lorsqu'aucune confusion n'est à craindre, nous noterons plus simplement $d$ cet opérateur, comme cela est l'usage.

Les définitions précédentes gardent un sens lorsque $U$ est un ouvert d'une variété compacte $C^\infty$. Prenons maintenant $V = \mathbf{T}^n$, et $m$ désigne toujours la partie entière de $\frac{n}{2}$. Notons $g_0$ la structure riemannienne standard sur $V$, et $\tau_0$ l'involution de Hodge associée. Soit $\delta_0 = \pm\tau_0 d\tau_0$ l'adjoint formel de $d$ sur $C^\infty(V, \Lambda_{\mathbf{C}}^m(T^*V))$ et $\Delta_0 = d\delta_0 + \delta_0 d$ le Laplacien. L'opérateur $\Delta_0$ est essentiellement autoadjoint sur $\Omega(V, g_0)$, son noyau est de dimension finie et ses éléments sont des fonctions $C^\infty$. Soit $P_0$ le projecteur de $L^2$ sur le noyau de $\Delta_0$ et notons $(P_0 + \Delta_0)^{-1}$ l'inverse de la fermeture de $P_0 + \Delta_0$ dans $L^2$. Si $\xi$ et $\eta$ sont des formes différentielles $C^\infty$ telles que $\partial\alpha + \partial\beta = n + 1$, l'égalité suivante est vraie:

$$\int_V \delta_0(P_0 + \Delta_0)^{-1}\xi \wedge \eta = \int_V \xi \wedge \delta_0(P_0 + \Delta_0)^{-1}\eta,$$

et montre que $\delta_0(P_0 + \Delta_0)^{-1}$ s'étend par dualité à l'espace de courants sur $V$.

Etant donné un espace de Hilbert $\mathcal{H}$, on note $\mathcal{L}(\mathcal{H})$ l'algèbre des opérateurs continus sur $\mathcal{H}$ et, $\mathcal{K}(\mathcal{H})$ l'idéal des opérateurs compacts. On note $\mathcal{L}^p(\mathcal{H})$ l'idéal des opérateurs compacts $T$ tels que $\mathrm{Tr}(\|T\|^p) < +\infty$, et pour $p > 1$, $\mathcal{L}^{p+}(\mathcal{H})$ l'idéal des opérateurs compacts $T$ tels que:

$$\sup_k \{(\lambda_1 + \cdots + \lambda_k)k^\alpha\} < +\infty$$

où $\alpha = p^{-1} - 1$ et $\lambda_1 \geq \cdots \geq \lambda_k \geq \cdots$ est la suite décroissante des valeurs propres de $T$.

Ces idéaux s'obtiennent par interpolation réelle à partir du couple $(\mathcal{K}(H), \mathcal{L}^1(\mathcal{H}))$ [4]. Plus généralement, étant donné un espace de Banach $E$, nous noterons $\mathcal{L}^{p+}(E)$ l'idéal d'opérateurs obtenus par interpolation réelle au point $(p^{-1}, 0)$ à partir du couple $(\mathcal{K}(E), \mathcal{L}^1(E))$ des idéaux d'opérateurs respectivement nucléaires et compacts sur $E$.

LEMME 1.5. *Soit $p, q$ et $N > n$ des réels positifs tels que $p^{-1} + q^{-1} = 1$ et $p^{-1} + N^{-1} = q^{-1}$ et soit $n_0 = nN(N-n)^{-1}$. Pour tout $\xi \in L^q(V, \Lambda_{\mathbf{C}}(T^*V))$, le courant $\delta_0(P_0 + \Delta_0)^{-1}\xi$ appartient à $L^p(V, \Lambda_{\mathbf{C}}(T^*V))$. L'application ainsi définie est continue et appartient à $\mathcal{L}^{n_0+}(L^q(V, \Lambda_{\mathbf{C}}(T^*V)), L^p(V, \Lambda_{\mathbf{C}}(T^*V)))$.*

*Démonstration.* Montrons que $\delta_0(P_0 + \Delta_0)^{-1}$ s'étend en un opérateur continu de l'espace $L^q(V, \Lambda_{\mathbf{C}}(T^*V))$ dans l'espace de Sobolev $L_1^q(V, \Lambda_{\mathbf{C}}(T^*V))$. Supposons d'abord que $n \geq 3$. Pour $0 < \alpha < n$, la transformée de Riesz est l'opérateur integro-singulier faible qui à $f \in \mathcal{S}(\mathbf{R}^n)$ associe:

$$I_\alpha(f) = \gamma(\alpha) \int_{\mathbf{R}^n} |x-y|^{\alpha-n} f(y) dy,$$

où $\gamma(\alpha) \in \mathbf{R}$ est une constante qui ne dépend que de $n$ et de $\alpha$ [10, Ch. 5]. Pour $1 < r < s < +\infty$, et $\frac{1}{s} + \frac{\alpha}{n} = \frac{1}{r}$, $I_\alpha$ s'étend en un opérateur continu de $L^s(\mathbf{R}^n) \to L^r(\mathbf{R}^n)$ et donc pour $\varphi \in C_c(\mathbf{R}^n)$, la multiplication par $\varphi$ est continue de $L^r(\mathbf{R}^n)$ dans $L^s(\mathbf{R}^n)$, et donc $\varphi I_\alpha$ est continu de $L^s$ dans $L^s$.

Soit $\varphi, \psi \in C_c^\infty(\mathbf{R}^n)$ telles que $\psi \equiv 1$ sur le support de $\varphi$. Montrons que $\psi I_2 \varphi$ s'étend par continuité en un opérateur de $L^q(\mathbf{R}^n)$ dans l'espace de Sobolev $L_2^q(\mathbf{R}^n)$. Soit $\Delta = -\sum \frac{\partial^2}{\partial x_i^2}$ le Laplacien scalaire et montrons que

(1.6) $$\Delta \psi I_2 \varphi = \varphi + k,$$

où $k \in C_c^\infty(\mathbf{R}^n)$. On a, pour $f \in \mathcal{S}(\mathbf{R}^n)$, $\Delta I_2(f) = f$ et:

$$\Delta \psi I_2 \varphi f = \varphi f + [\Delta, \psi] I_2 \varphi f = \varphi f + \Delta(\psi) I_2 \varphi f + \sum \frac{\partial \psi}{\partial x_i} \frac{\partial}{\partial x_i}(I_2 \varphi f).$$

Remarquons que le noyau distribution de $I_2$ est $C^\infty$ en dehors de la diagonale. Comme $\text{supp}(\Delta(\psi)) \cap \text{supp}(\varphi) = \emptyset$, $[\Delta, \psi] I_2 \varphi$ est un noyau $C^\infty$, et pour la même raison $\frac{\partial \psi}{\partial x_i} I_2 \varphi$ est un noyau $C^\infty$ et donc $\frac{\partial \psi}{\partial x_i} \frac{\partial}{\partial x_i}(I_2 \varphi)$ l'est aussi.

Soit $g \in L^q$ telle que $g = \psi I_2 \varphi f$ pour $f \in L^q$. On a alors $\Delta g = \varphi f + kf$, ce qui montre que $g \in L_2^q$ et la formule (1.6) montre que l'application ainsi définie est continue de $L^q$ dans $L_2^q$ [10, Cor., p. 77].

Soit $U_1, U_2$ deux ouverts de carte de $V$ tels que $U_1 \cup U_2 = V$ et $\varphi_1, \varphi_2$ une partition de l'unité subordonnée à ce recouvrement et $\psi_i \in C_c^\infty(U_i)$ telles que $\psi_i \varphi_i = \varphi_i$ et soit l'opérateur $T_i = \psi_i I_2 \varphi_i$ agissant sur $C^\infty(U_i)$ et $T = (T_1 + T_2) \otimes 1$ agissant sur $C^\infty(V, \Lambda_{\mathbf{C}}(T^*V))$. La formule (1.6) montre que $\Delta_0 T = 1 + S$ où $S$ est un noyau régularisant et donc $T - (P_0 + \Delta_0)^{-1}$ est aussi régularisant. En particulier, $(P_0+\Delta_0)^{-1}$ s'étend en un opérateur continu de $L^q(V, \Lambda_{\mathbf{C}}(T^*V))$ dans l'espace de Sobolev $L_2^q(V, \Lambda_{\mathbf{C}}(T^*V))$ et donc $\delta_0(P_0+\Delta_0)^{-1}$ s'étend en un opérateur continu de $L^q(V, \Lambda_{\mathbf{C}}(T^*V))$ dans $L_1^q$. Cela montre l'assertion pour $n \geq 3$; pour $n < 3$, on refait le même raisonnement à partir des noyaux des solutions fondamentales de $\Delta$, ou bien se ramène au cas précédent en prenant $V \times \mathbf{T^3}$.

Par le lemme de Sobolev, on une inclusion continue $t : L_1^q(V, \Lambda(T^*V)) \to L^p(V, \Lambda(T^*V))$, puisque $p^{-1} > q^{-1} - n^{-1}$. Cette inclusion $t$ est compacte, par le lemme de Rellich. Pour montrer la dernière assertion, il suffit alors de montrer que $t \in \mathcal{L}^{n_0+}$. La transformation de Fourier induit un isomorphisme

$$\mathcal{F} : C^\infty(V, \Lambda_{\mathbf{C}}(T^*V)) \to \mathcal{S}(\mathbf{Z}^n, \Lambda_{\mathbf{C}}(\mathbf{R}^{n*}))$$

et par le théorème de Hausdorff-Young, $\mathcal{F}$ s'étend par continuité en un opérateur borné de $L^q(V, \Lambda_{\mathbf{C}}(T^*V))$ dans $l^r(\mathbf{Z}^n, \Lambda_{\mathbf{C}}(\mathbf{R}^{n*}))$ et $\mathcal{F}^{-1}$ s'étend par continuité en un opérateur borné de $l^s(\mathbf{Z}^n, \Lambda_{\mathbf{C}}(\mathbf{R}^{n*}))$ dans $L^p(V, \Lambda_{\mathbf{C}}(T^*V))$, avec:

$$r = \frac{q}{q-1} \qquad \text{et} \qquad s = \frac{p}{p-1}.$$

La restriction de $\mathcal{F}$ à $L_1^q(V, \Lambda_{\mathbf{C}}(T^*V))$ a pour image le sous-espace $l_1^r(\mathbf{Z}^n, \Lambda_{\mathbf{C}}(\mathbf{R}^{n*}))$ de $l^r$ formé par les suites $\xi$ telles que $\|\xi\|_{1,r} = (\sum_k (1 + \|k\|)^r |\xi(k)|^r)^{\frac{1}{r}} < +\infty$. Pour tout $K > 0$, soit $t_K : l_1^r \to l^s$ donné par:

$$t_K \xi(k) = \begin{cases} \xi(k) & \text{si} \quad \|k\| \geq K \\ 0 & \text{si} \quad \|k\| < K. \end{cases}$$

On a $\|t - t_K\|_1 = \mathcal{O}(1+K)^{n-1}$, et l'inégalité de Hölder montre que $\|t_K\|_\infty \leq \{\sum_{\|k\| \geq K} (1 + \|k\|)^{-N}\}^{\frac{1}{N}} = \mathcal{O}((1+K)^{\frac{n}{N}-1})$, $K \to +\infty$. Avec les notations de [3], nous avons alors l'estimé pour $\lambda \to +\infty$:

$$K(\lambda, t) = \mathcal{O}(\lambda^{\frac{1}{N} - \frac{1}{n}})$$

ce qui montre que $t \in \mathcal{L}^{n_0+}(l^r, l^s)$. $\square$

Pour $g \in \mathcal{R}(V)$, soit $n(g)$ le réel:

$$n(g) = \frac{n p_m q_{m+1}}{p_m q_{m+1} - n(p_m - q_{m+1})}.$$

Par le (2) de la définition 1.1, on a l'inégalité $n(g) \geq 2$ avec égalité si et seulement si $p_m = q_{m+1}$, c'est à dire si $\Omega^m(V,g)$ est égal en tant qu'espace de Banach à l'espace de Hilbert $L^2(V, \Lambda_{\mathbf{C}}^k(T^*V))$.

Cela se produit si $g = g_0$ est la structure riemannienne plate standard et plus généralement si $g$ est une structure $C^\infty$. Si $n$ est pair cela lieu aussi pour les structures quasi-conformes [5]. En général, ni l'inclusion $\Omega^m(V,g) \subset L^2(V, \Lambda_{\mathbf{C}}^m(T^*V))$, ni l'inclusion opposée n'ont lieu.

Pour la proposition suivante, nous prenons toujours $V = \mathbf{T}^n$ et fixons $g \in \mathcal{R}(V)$. Nous avons donc un opérateur fermé et densément défini $d : \Omega^m(V,g) \to \Omega^{m+1}(V,g)$.

PROPOSITION 1.7. *Si $n$ est impair, l'opérateur $D = \tau d$ agissant sur $\Omega^m(V,g)$ est autoadjoint et $D(1+D^2)^{-1}$ appartient à l'idéal $\mathcal{L}^{n(g)+}(\Omega^m(V,g))$.*

*Si $n$ est pair, l'opérateur $D = d + d^* : \Omega^m \to \Omega^{m-1} \oplus \Omega^{m+1}$ anticommute à $\tau$ et $(1+|D|)^{-1} \in \mathcal{L}^{n(g)+}(\Omega^m(V,g))$.*

*Démonstration.* Prenons d'abord $n$ impair. Par le lemme précédent, l'opérateur $t$ déterminé par l'égalité $t\omega = \delta_0(P_0 + \Delta_0)^{-1}\omega$ appartient à $\mathcal{L}^{n(g)+}(\Omega^{m+1}(V,g), \Omega^m(V,g))$.

Soit $P$ le projecteur orthogonal de $\Omega^m(V,g)$ sur le support de $d$, c'est à dire sur le sous-espace orthogonal à $\ker d$, et soit $Q$ le projecteur orthogonal de $\Omega^{m+1}(V,g)$ sur la fermeture de $\operatorname{im} d$. Comme la relation $dtd\omega = d\omega$ est vraie au sens des distributions, elle reste vraie pour $\omega \in \operatorname{dom} d$. L'opérateur $d^{-1}$, inverse de $d$ sur son image est alors donné par $d^{-1} = PtQ$, et est donc élément de $\mathcal{L}^{n(g)+}(\Omega^m(V,g))$.

Montrons que $\tau d$ est autoadjoint. Soit $d_0$ la fermeture de la restriction de $d$ aux formes différentielles $C^\infty$ sur $V$. Le graphe du transposé de $d_0$ est, dans le cas $n$ impair, l'ensemble des $(\beta, \gamma) \in \Omega^m(V,G) \oplus \Omega^{m+1}(V,g)$ tels que pour toute forme $\alpha$ de classe $C^\infty$, on ait:

$$\int_V d\alpha \wedge \beta = (-1)^{m-1} \int_V \alpha \wedge \gamma.$$

Comme la différentielle extérieure de de Rham usuelle pour une forme $C^\infty$ est égale à sa différentielle extérieure au sens des distributions, cela entraîne $\beta \in \operatorname{dom} d$ et $\gamma = d\beta$, et donc ${}^t d_0 = d$. L'opérateur adjoint de $\tau d$ est alors:

$$(\tau d)^* = \tau^t(\tau d)\tau = \tau^t d^t \tau \tau = \tau d_0.$$

Il suffit donc d'établir que $d = d_0$. Soit $d_1$ le fermeture de la restriction de $d$ au sous-domaine $\operatorname{dom} d \cap L^{p_m}$. Comme $\Omega^{m+1}(V,g) \subset L^{q_{m+1}}$, $d_1$ est un opérateur symétrique, i.e. pour $\alpha, \beta \in \operatorname{dom} d_1$, on a:

$$\int_V \alpha \wedge d\beta = (-1)^{m-1} \int_V d\alpha \wedge \beta.$$

Cette égalité est satisfaite puisque:

$$\frac{1}{p_m} + \frac{1}{q_{m+1}} \leq 1.$$

Si $\xi \in \operatorname{im} d$, alors $t\xi \in L^{p_m}$ et $dt\xi = \xi$, ce qui montre que $\operatorname{im} d_1 = \operatorname{im} d$. En particulier $d_1 = {}^t d_1$: en effet, $\ker d_1 = \ker {}^t d_1$ et donc pour tout $\alpha \in \operatorname{dom} {}^t d_1$, il existe $\beta \in \operatorname{dom} d_1$ tel que $d_1 \beta = {}^t d_1 \alpha$ et $\alpha - \beta \in \ker {}^t d_1$, ce qui montre que $\alpha \in \ker {}^t d_1 + \operatorname{dom} d_1 \subset \operatorname{dom} d_1$. Il suffit alors de montrer que $\ker d = \ker d_1$.

En approchant la norme hilbertienne sur $\Omega^m$ par une suite croissante $\|\cdot\|_r$, $r \in \mathbf{N}$ de normes hilbertiennes fibrées $L^\infty$, on obtient une famille continue croissante d'espaces de Hilbert $\mathcal{K}_r \subset \mathcal{K}_{r+1} \subset \cdots \Omega^m(V,g)$ pour laquelle la suite des opérateurs $d_1$ est continue. Mais, dans $\mathcal{K}_r$, $\ker d$ est dans l'adhérence de $\ker d_1$: cela résulte d'un raisonnement semblable appliqué au dual de $\mathcal{K}_r$, la famille croissante étant ici topologiquement standard.

Le cas où $n = 2m$ est un petit peu différent. Comme précédemment, soit $t_k \in \mathcal{L}^{n(g)+}(\Omega^k(V,g), \Omega^{k-1}(V,g))$, pour $k = m$ et $k = m+1$, l'opérateur déterminé par l'égalité $t_k\omega = \delta_0(P_0 + \Delta_0)^{-1}\omega$. En notant $P$ le projecteur orthogonal de $\Omega^m(V,g)$ sur le support de $d : \Omega^m(V,g) \to \Omega^{m+1}(V,g)$, et $Q$ le projecteur orthogonal de $\Omega^m(V,g)$ sur $\mathrm{im}\, d$, nous obtenons de même que les opérateurs $(P + d^*d)^{-\frac{1}{2}}$ et $(Q + dd^*)^{-\frac{1}{2}}$ appartiennent à $\mathcal{L}^{n(g)+}(\Omega^m(V,g))$ et nous avons ${}^td = d$, et donc $d^* = -\tau d\tau$.

Il ne reste plus qu'à montrer que le projecteur $R = 1 - P - Q$ est de rang fini. Soit $\xi, \eta \in C^\infty$, et $T(\xi,\eta)\zeta = (\int_V \zeta \wedge \tau_0\bar\eta)\xi$, c'est à dire l'opérateur de rang un associé pour la structure riemannienne standard. Cet opérateur s'étend par continuité en un opérateur continu de rang un de $\Omega^m(V,g)$. En effet, comme $\tau_0\eta \in \Omega^m(V,g)$, on a $\int_V \zeta \wedge \tau_0\bar\eta = \int_V \tau\zeta \wedge \tau\tau_0\bar\eta$ et donc $\|T(\xi,\eta)\| \leq \|\xi\|\|\tau_0\eta\|$. Nous avons alors:
$$dt_m\xi + t_{m+1}d\xi = \xi + a\xi$$
où $a$ est une combinaison linéaire finie de $T(\xi,\eta)$, et donc détermine un opérateur de rang fini de $\Omega^m(V,g)$. Cette égalité est encore vraie au sens des distributions et s'applique aux courants. Soit $\alpha \in \mathrm{im}\, R$ ; comme $\Omega^m(V,g) \subset L^{q_m}$, nous avons alors $dt_m\alpha = \alpha + a\alpha$, ce qui montre que $t_m\alpha \in \mathrm{dom}\, d$, et qu'en particulier on peut appliquer $R$ à $dt_m\xi$, ce qui donne $0$, et donc $\mathrm{im}\, R$ est un sous-espace de $\ker(1 + Ra)$ qui est de dimension finie. □

## 2. Variétés topologiques ayant une structure $L^p$

Soit $\theta : U \to W$ un homéomorphisme entre ouverts de $\mathbf{R}^n$ dérivable presque partout et dont la dérivée $\theta'$ appartient $L^p(U, \mathrm{End}(\mathbf{R}^n))$ pour un $p \geq 1$. Nous dirons que $\theta$ est dérivable d'ordre $p$. Un tel homéomorphisme préserve la classe de la mesure de Lebesgue [9] et la dérivée agit sur les sections mesurables du fibré tangent.

*Définition* 2.1. *Une variété topologique compacte* est dite dérivable à l'ordre $p$ si les changements de carte associés à son atlas sont dérivables d'ordre $p$.

Il est possible alors de définir des structures riemanniennes sur $V$ et l'espace de Hilbert des formes différentielles mesurables sur $V$ de degré $k$ et de carré intégrable pour une structure $g$ sera noté $\Omega^k(V,g)$.

Etant donné un recouvrement $\mathcal{O} = (\mathcal{O}_i)_{i \in I}$ par des ouverts de carte $\theta_i : \mathcal{O}_i \to U_i \subset \mathbf{R}^n$, nous noterons $\mathcal{R}(\mathcal{O})$ les structures riemanniennes $g$ telles que, en notant $g_i$ l'image sur $U_i$ de la restriction de $g$ à $\mathcal{O}_i$, on ait $g_i \in \mathcal{R}(U_i)$, et on posera

$$n(g) = \sup_i n(g_i).$$

Soit $\mathcal{B}(\mathcal{O}) \subset C(V)$ la sous-algèbre dense engendrée par les sous-algèbres $C_c^\infty(U_i)$, pour $i \in I$, c'est à dire engendrée par les éléments de la forme $f = f_1 + \cdots + f_N$ où $f_i \in C_c^\infty(U_i)$. Si $g \in \mathcal{R}(\mathcal{O})$, la relation (3) de la définition 1.1 nous montre que le sous-espace engendré par la réunion des $C_0(U_i, \Lambda_{\mathbf{C}}^k(T^*U_i))$ est un sous-espace dense de $\Omega^k(V, g)$. L'espace $\mathcal{R}(\mathcal{O})$ devient alors un espace métrisable pour la famille de semi-distances déterminées par les fonctions positives:

$$g \to \|\omega\|_g^2,$$
$$g \to p_k(g_i),$$
$$g \to q_k(g_i),$$

où $\omega$ parcourt la réunion des $C_0(U_i, \Lambda_{\mathbf{C}}^k(T^*U_i))$ et $p_k(g_i)$ (resp. $q_k(g_i)$), pour $k = m, m+1$ est le plus petit (resp. plus grand) réel pour lesquels les inclusions de la définition 1.1 sont vérifiées sur $U_i$. Pour cette topologie,

PROPOSITION 2.2. *L'espace topologique métrisable $\mathcal{R}(\mathcal{O})$ est connexe par arcs.*

*Démonstration.* Cette assertion est une conséquence de l'interpolation complexe. Soit $U$ un ouvert relativement compact de $\mathbf{R}^n$, et $g_0, g_1$ deux éléments de $\mathcal{R}(U)$, $p_k^0, p_k^1, q_k^0, q_k^1$ les réels associés vérifiant la condition (1) de la définition 1.1. Il y a une inclusion continue de $\Omega(U, g_i)$ ($i = 0, 1$) dans $L^q(U, \Lambda_{\mathbf{C}}(T^*(U)))$ avec $q = \min(q_0, q_1)$, et nous pouvons appliquer l'interpolation complexe au couple d'espaces de Hilbert $(\Omega^*(U, g_0), \Omega^*(U, g_1))$. Pour $X \in \mathbf{R}^n$, notons $\langle X, X \rangle$ le produit scalaire euclidien standard et, pour $i = 0$ et $i = 1$, soit $x \to A_i(x)$ le champ mesurable sur $U$ de matrices positives tel que pour tout $X \in \mathbf{R}^n$, on ait l'égalité presque partout pour la mesure de Lebesgue $g_i(x)(X, X) = \langle A_i(x)X, X \rangle$. Posons, pour $t \in [0, 1]$:

$$A_t = A_0^{\frac{1}{2}}(A_0^{-\frac{1}{2}} A_1 A_0^{-\frac{1}{2}})^t A_0^{\frac{1}{2}}.$$

Nous avons alors une identification canonique de $(\Omega^k(U, g_0), \Omega^k(U, g_1))_t$ avec $\Omega(U, g_t)$ où $g_t(x)(X, X) = \langle A_t(x)X, X \rangle$. Pour $\omega \in C^\infty(U, \Lambda_{\mathbf{C}}^k(\mathbf{R}^{n,*}))$, on a en effet pour tout $t \in [0, 1]$:

$$\|\omega\|_{g_t}^2 = \int_U \langle \lambda^k(A_t^{-1}(x))\omega, \omega \rangle \sqrt{\det(A_t(x))} dx.$$

En appliquant la méthode de [3, Ch. 5], on a canoniquement que, pour $\omega$ de degré $k$, la norme sur $(\Omega^*(U,g_0),\Omega^*(U,g_1))_t$ est donnée par:

$$\int_U \langle B_t(x)\omega(x),\omega(x)\rangle dx$$

où $B_t = \lambda^k(A_0)^{-\frac{1}{2}}(\lambda^k(A_0)^{\frac{1}{2}}\lambda^k(A_1^{-1})\lambda^k(A_0)^{\frac{1}{2}})^t\lambda^k(A_0)^{-\frac{1}{2}}\det(A_0(x)^{1-t}A_1(x)^t)^{\frac{1}{2}}$, c'est à dire $B_t = \lambda^k(A_t^{-1})\sqrt{\det(A_t(x))}$.

Les relations de la définition 1.1 sont satisfaites pour $g(t)$ avec ([3, Ch. 4]):

$$\frac{1}{p_k(t)} = \frac{t}{p_k^0} + \frac{1-t}{p_k^1},$$
$$\frac{1}{q_k(t)} = \frac{t}{q_k^0} + \frac{1-t}{q_k^1}.$$

L'application $t \to g(t)$ est bien continue. □

*Remarque* 2.3. L'opérateur $A_t = A_0^{\frac{1}{2}}(A_0^{-\frac{1}{2}}A_1 A_0^{-\frac{1}{2}})^t A_0^{\frac{1}{2}}$, pour $t \in [0,1]$, est un cas particulier de la notion de *moyenne d'opérateurs linéaires positifs* introduite par T. Ando et F. Kubo [1].

*Remarque* 2.4. Soit $f : Z \to \mathcal{R}(\mathcal{O})$ une application continue: sur le champ $\mathcal{E} = (\Omega(V,f(t))_{t\in Z}$, il existe une unique structure de champ continu engendrée par les sections $\omega$, $\omega \in \sum_i C_0(U_i, \Lambda^k_{\mathbf{C}}(T^*U_i))$. On a donc sur $\mathcal{E}$ une structure de $C_0(Z)$-module hilbertien canoniquement associé à $f$.

Fixons nous un tel recouvrement $\mathcal{O}$ de $V$ et un élément $g \in \mathcal{R}(\mathcal{O})$ et, comme précédemment, soit $m$ la partie entière de $\frac{n}{2}$. Commençons par définir $d : \Omega^k(V,g) \to \Omega^{k+1}(V,g)$. Soit $\omega \in \Omega^k$ et $\omega_i$ l'image dans $U_i$ de sa restriction à $\mathcal{O}_i$. Nous dirons que $\omega \in \mathrm{dom}\, d$ si pour tout $i$, $\omega_i \in \mathrm{dom}\, d_{U_i} \subset \Omega^{k+1}(U_i,g_i)$ et on définit alors $d\omega$ par la condition que la restriction de $d\omega$ à $U_i$ soit égal à $d_{U_i}\omega_i$. Le lemme suivant nous montre que cette définition a un sens et que $d$ a un domaine dense.

LEMME 2.5.　　*Soit $U$ et $W$ deux ouverts de $\mathbf{R}^n$ et $\theta : U \to W$ un homéomorphisme dérivable d'ordre $p \geq m+1$, et $g_1 \in \mathcal{R}(U)$, $g_2 \in \mathcal{R}(W)$ des structures riemanniennes telles $\theta^*(g_2) = g_1$, et $\alpha \in \Omega^k(W,g_2)$ telle que $d_W\alpha \in \Omega^{k+1}$. Les propriétés suivantes sont vraies pour $k = m-1, m$ si $n$ est pair et pour $k = m$ si $n$ est impair*:

1) Il existe une suite $\beta_n \in C_c^\infty(\mathbf{R}^n, \Lambda^k(T^*\mathbf{R}^{n,*}))$ telle que si $\alpha_n$ est la restriction de $\beta_n$ à $W$, on a $\lim \alpha_n = \alpha$ et $\lim d_W\alpha_n = d_W\alpha$.
2) $\theta^*(\alpha) \in \mathrm{dom}\, d_U \subset \Omega(U,g_1)$ et $d_U\theta^*(\alpha) = \theta^*(d_W\alpha)$.

*Démonstration*. Supposons $n$ impair. Si le support de $\alpha$ est compact dans $W$, le 1) résulte alors du cas du tore. Le cas général s'en déduit alors par transposition. Pour montrer 2), il suffit alors de supposer que $\alpha$ est la restriction à $W$ de $\beta \in C_c^\infty(\mathbf{R}^n, \Lambda^k(T^*\mathbf{R}^{n,*}))$ puisque $d$ est fermé. Comme $\theta' \in L^p$ avec $p \geq m+1$, les formes différentielles $\theta^*(\alpha)$ et $\theta^*(d\alpha)$ sont dans $L^1$, et sont donc des courants. Soit alors $\theta_\varepsilon : U \to \mathbf{R}^n$ une suite d'applications $C^\infty$ (non nécessairement bijectives) convergeant vers $\theta$ uniformément et dont les dérivées $\theta'_\varepsilon$ convergent vers $\theta'$ dans $L^p$. Pour $\zeta \in C_c^\infty(U, \Lambda^m(T^*U))$, il est vrai que:

$$\int_U \theta^*(\alpha) \wedge d\zeta = \lim \int_U \theta_\varepsilon^*(\beta) \wedge d\zeta = \lim (-1)^{m-1} \int_U \theta_\varepsilon^*(d\beta) \wedge \zeta$$
$$= (-1)^{m-1} \int_U \theta^*(d\alpha) \wedge \zeta. \qquad \square$$

2.6. Pour le prochain théorème, nous nous plaçons dans les conditions suivantes: $V$ est une variété compacte, orientée, de dimension $n$, admettant une structure dérivable d'ordre $m+1$ et telle que pour tout recouvrement fini par des ouverts de carte $\theta_i : \mathcal{O}_i \to U_i$, l'espace $\mathcal{R}(\mathcal{O})$ est non vide. Pour une telle variété, nous pouvons construire analytiquement un opérateur de signature unique à homotopie près:

THÉORÈME 2.7. *Avec les hypothèses précédentes, soit un recouvrement $\mathcal{O}$ de $V$ et $g \in \mathcal{R}(\mathcal{O})$.*

*Si $n$ est impair, l'opérateur $D = \tau d$ sur $\Omega^m(V, g)$ est densément défini, fermé, autoadjoint et a une résolvante sur son support dans $\mathcal{L}^{n(g)+}$.*

*Si $n$ est pair, l'opérateur $D = d + d^*$ de $\Omega^m(V, g)$ dans $\Omega^{m-1}(V, g) \oplus \Omega^{m+1}(V, g)$ est densément défini, fermé, anticommute à $\tau$ et a une résolvante dans $\mathcal{L}^{n(g)+}$.*

*Quelle que soit la parité de $n$, pour tout $f \in \mathcal{B}(\mathcal{O})$, le commutateur $[D, f]$ est densément défini et borné.*

*L'opérateur $D$ obtenu est unique à homotopie non bornée près.*

La dernière assertion du théorème veut dire, par exemple dans le cas impair, que si $\mathcal{O}_1$ est un autre recouvrement par des ouverts de carte de $V$ et $g_1 \in \mathcal{R}(\mathcal{O}_1)$, il existe alors un champ continu $\mathcal{E}$ d'espaces de Hilbert sur $[0,1]$ tel que $\mathcal{E}_0 = \Omega^m(V, g)$, $\mathcal{E}_1 = \Omega^m(V, g_1)$, et une famille continue d'opérateurs autoadjoints $\mathcal{D}$ sur $\mathcal{E}$ à résolvante dans l'idéal des opérateurs compacts du module hilbertien, telle que $\mathcal{D}_0 = D$, $\mathcal{D}_1 = D_1$ et une sous-algèbre involutive dense de $C(V \times [0,1])$ dont les éléments commutent à $\mathcal{D}$ modulo les bornés. Dans le cas pair, on a une description analogue en remplaçant $\Omega^m(V, g)$ par $\Omega^m(V, g) \oplus \Omega^{m-1}(V, g) \oplus \Omega^{m+1}(V, g)$.

*Démonstration.* Soit $\varphi_i \in \mathcal{B}(\mathcal{O})$ une partition de l'unité affiliée au recouvrement. Le domaine de $d$ contient le sous-espace de $\Omega^m$ engendré par les sous-espaces $C_c^\infty(U_i, \Lambda^m)$ qui est dense. Soit $\omega \in \mathrm{dom}\, d$ et $\omega_i = \varphi_i \omega$. Par le lemme précédent, il existe une suite $\alpha_k \in C_c^\infty(U_i, \Lambda^m)$ qui converge vers $\omega_i$ et telle que $d\alpha_k$ converge vers $d\omega_i$: cela montre que, si $n$ est impair, $\tau d$ est autoadjoint et si $n$ est pair, que l'adjoint de $d$ est $-\tau d \tau$. Supposons maintenant que $n$ est impair et montrons que $\mathrm{im}\, d$ est fermée. Soit $d_i$ la fermeture de l'opérateur $d$ sur le domaine essentiel $C_c^\infty(\mathcal{O}_i, \Lambda_{\mathbf{C}}^m(T^*\mathcal{O}_i))$ et ${}^t d_i$ son transposé sur $\Omega^m(\mathcal{O}_i)$. Par ce qui précède, on sait que $\mathrm{im}\, d_i$ est fermée et donc par dualité, $\mathrm{im}\, {}^t d_i$ l'est aussi. Soit $\omega_k \in \mathrm{dom}\, d \subset \Omega^m(V, g)$ une suite telle que $d\omega_k$ converge. On sait alors qu'il existe $\xi_i \in \mathrm{dom}\, {}^t d_i$ tel que ${}^t d_i \xi_i = \lim(d\omega_{|\mathcal{O}_i})$. Posons alors $\xi = \sum_i \varphi_i \xi_i$; on a bien $\lim d\omega_k = d\xi$, et donc $\mathrm{im}\, d$ est fermée, et l'inverse $d^{-1}$ de $d$ sur son support est continu.

Soit $P$ le projecteur orthogonal de $\Omega^m(V, g)$ sur le support de $d$, et $Q$ celui de $\Omega^{m+1}$ sur $\mathrm{im}\, d$. Pour chaque $i$, soit $\mathcal{O}_i \to U_i \subset \mathbf{T}^n$ une immersion et fixons une extension $h_i$ sur $\mathbf{T}^n$ de $g_i$ sur $U_i$. Nous identifions alors les formes différentielles à support dans $\mathcal{O}_i$ et celles sur $\mathbf{T}^n$ à support dans $U_i$. Soit $t^i$ l'opérateur continu de $\Omega^{m+1}(\mathbf{T}^n, h_i) \to \Omega^m(\mathbf{T}^n, h_i)$ défini dans le lemme 1.4, et tel que, pour $\xi$ de classe $C^\infty$:

$$t^i \xi = \delta_0 (P_0 + \Delta_0)^{-1} \xi.$$

Soit $\psi_i \equiv 1$ sur $\mathrm{supp}(\varphi_i)$ et à support dans $U_i$, et $[d, \psi_i]$ borné et notons $d_i$ la différentielle extérieure sur $\Omega(\mathbf{T}^n, h_i)$. Soit alors $l_i : \Omega^{m+1}(V, g) \to \Omega^{m+1}(\mathbf{T}^n, h_i)$ donné par $l_i(\alpha) = d_i \psi_i d^{-1} Q \alpha$ et posons $T\alpha = \sum \varphi_i t_i l_i(\alpha)$. Calculons $dT$; on a:

$$dT\alpha = \sum_i \varphi_i l_i(\alpha) + \sum_i [d, \varphi_i] t_i l_i(\alpha).$$

Comme $\varphi_i[d, \psi_i] = 0$, et $\varphi_i \psi_i = \varphi_i$, on a $dT\alpha = Q + bQ$ où $b$ est donné par le dernier terme de l'égalité ci-dessus et appartient à $\mathcal{L}^{n(g)+}(\Omega^{m+1}, \Omega^m)$. Nous avons alors $PT - d^{-1} \in \mathcal{L}^{n(g)+}$ et donc $d^{-1} \in \mathcal{L}^{n(g)+}$.

Dans le cas où $n$ est pair, nous obtenons par un raisonnement analogue que l'image de $d$ est fermée et que $d = d_0$ et $d^* = -\tau d \tau$ et donc $d + d^*$ anticommute à $\tau$. On construit les opérateurs $T^k \in \mathcal{L}(\Omega^{k+1}, \Omega^k)$, $T^k \alpha = \sum \varphi_i t_i^k l_i^k(\alpha)$, où $t_i^k, l_i^k$ sont construits comme précédemment. Nous obtenons que $dT^{m-1} + T^m d = 1 + a$ avec $a \in \mathcal{L}^{n(g)+}$, que $d^{-1}, d^{*-1} \in \mathcal{L}^{n(g)+}$ et finalement que $D$ est à résolvante dans $\mathcal{L}^{n(g)+}$. Notons que cette propriété peut aussi se démontrer comme dans [9] à l'aide de l'inclusion canonique du graphe de $D$ dans $\Omega(V, g)$.

Pour la dernière assertion, soit $\mathcal{O}_2$ le recouvrement le moins fin qui contient $\mathcal{O}$ et $\mathcal{O}_1$ et soit $h_2 \in \mathcal{R}(\mathcal{O}_2)$. Il suffit alors de démontrer l'assertion pour $h_2$. Soit $g(t)$ un chemin continu dans $\mathcal{R}(\mathcal{O})$ avec $g(0) = h_2$ et $g(1) = g$ obtenu par interpolation complexe comme dans la preuve de la 2.2. Comme $\tau d$ est essentiellement autoadjoint sur le sous-espace engendré par $C_c^\infty(U_i, \Lambda^m)$, la proposition 2.9 de [8] et les théorèmes généraux d'interpolation montrent que, avec des notations évidentes, la famille d'opérateurs $\mathcal{D} = (D_t)_{t\in[0,1]}$ agissant sur le C([0,1])-module $(\Omega^*(V, g(t)))_{t\in[0,1]}$ est autoadjointe et à résolvante dans l'idéal des opérateurs compacts du module.

Soit enfin $W_j$, pour $1 \leq j \leq K$ les ouverts du recouvrement $\mathcal{O}_2$. La sous-algèbre de $C(V \times [0,1])$ engendrée par la réunion des sous-algèbres

$$C_c^\infty(U_i \times [0,1]) \quad \text{et} \quad C_c^\infty(W_j \times ]0,1])$$

est dense, est constituée d'éléments qui commutent à $\mathcal{D}$. $\square$

## 3. Exemples

La façon la plus directe de construire une structure riemannienne consiste à choisir un recouvrement par des ouverts de cartes $\theta_i : \mathcal{O}_i \to U_i \subset \mathbf{R}^n$, $i \in \{1, \cdots, k\}$. Pour $i, j \in \{1, \cdots, k\}$, soit alors $\psi_{i,j} : U_i \to M_n(\mathbf{R})$ la fonction obtenue en prolongeant la dérivée du changement de carte, par 0 en dehors de son domaine de définition naturel, c'est à dire:

$$\psi_{i,j}(x) = \begin{cases} d_x(\theta_j \circ \theta_i^{-1}) & \text{si } x \in \theta_i(\mathcal{O}_i \cap \mathcal{O}_j) \\ 0 & \text{sinon} \end{cases}$$

et posons, pour $X \in \mathbf{R}^n$ et $x \in \mathcal{O}_i$:

$$g_i(x)(X) = \|X\|^2 + \sum_{j \neq i} \|\psi_{i,j} X\|^2.$$

Comme $(\theta_j \circ \theta_i^{-1})^* g_j = g_i$, nous avons bien une structure riemannienne sur $V$ et telle que $g_i \geq 1$. Nous appellerons cette structure riemannienne *métrique spécifique du recouvrement*. Voici deux exemples de variétés pour lesquelles à chaque recouvrement ouvert $\mathcal{O}$, la métrique spécifique appartient à $\mathcal{R}(\mathcal{O})$, ce qui nous place dans les conditions du théorème 2.7.

3.1. *Variétés quasi-conformes.* Soit $\theta : U \to W$ un homéomorphism quasi-conforme entre deux ouverts de $\mathbf{R}^n$. Il a été montré par W. Gehring [6] que l'application $\theta$ est dérivable et qu'il existe $p > n$ tel que $\theta' \in L^p(U, \text{End}(\mathbf{R}^n))$. Par hypothèse, il existe alors $\varphi \in L^p(U)$ et $a, b > 0$ tels que pour presque tout $x \in U$, $\text{Sp}^t \theta'(x)\theta'(x) \subset [\varphi^2 a, \varphi^2 b]$.

Soit $V$ une variété quasi-conforme compacte et $(\mathcal{O}_i)_{i \in I}$ un recouvrement fini par des ouverts de carte $\theta_i : \mathcal{O}_i \to U_i$ et $g = (g_i)_{i \in I}$ la métrique spécifique de ce recouvrement. Soit $i$ fixé, il existe $p > n$, une structure riemannienne $g^0$ de classe $C^\infty$, et une fonction $\varphi \geq 1$, $\varphi \in L^p(U_i)$, tels que $g_i$ soit équivalente à $\varphi g^0$.

En tant qu'espaces de Banach, on a alors

$$\Omega^k(U_i, g) = L^2(U_i, \Lambda_\mathbf{C}^k(\mathbf{R}^n), \varphi^{n-2k} dx).$$

Les relations de la définition 1.1 sont donc satisfaites avec, pour $n = 2m + 1$:

$$p_m = \frac{2p}{(p-1)}, \qquad q_m = p_{m+1} = 2, \qquad q_{m+1} = \frac{2p}{(p+1)}$$

et pour $n = 2m$:

$$p_{m-1} = \frac{2p}{p-2}, \quad p_m = q_m = q_{m-1} = p_{m+1} = 2, \quad q_{m+1} = \frac{2p}{p+2}.$$

Vérifions-le si $n = 2m+1$; pour toute fonction mesurable $\psi$, l'inégalité de Hölder nous donne:

$$\left\{ \int_{U_i} |\psi(x)|^2 \varphi(x) dx \right\}^{\frac{1}{2}} \leq \left\{ \int_{U_i} |\psi(x)|^{p_m} dx \right\}^{\frac{1}{p_m}} \left\{ \int_{U_i} \varphi(x)^p dx \right\}^{\frac{1}{2p}}$$

et comme $\varphi \geq 1$, nous avons bien les inclusions:

$$L^{p_m}(U_i, \Lambda_{\mathbf{C}}^m(\mathbf{R}^n), dx) \subset L^2(U_i, \Lambda_{\mathbf{C}}^m(\mathbf{R}^n), \varphi dx) \subset L^2(U_i, \Lambda_{\mathbf{C}}^m(\mathbf{R}^n), dx).$$

Comme $p > n$, on a enfin:

$$\frac{1}{p_m} + \frac{1}{n} = \frac{1}{2}(1 - \frac{1}{p} + \frac{2}{n}) > \frac{1}{q_{m+1}}.$$

3.2. *Variétés dérivables à l'ordre p.* Si $V$ est une variété topologique dont les changements de carte de son atlas sont des homéomorphismes dérivables $\theta : U \to W$ avec $\theta' \in L^p$. Les conditions énoncées du théorème 2.7 sont satisfaites dès que:

$$p > \frac{n(n+1)}{2}.$$

En effet, avec les notations précédentes, soit $(\mathcal{O}_i)_{i \in I}$ un recouvrement et $(g_i)_{i \in I}$ sa métrique spécifique. Soit $A : U_i \to \mathcal{L}(\mathbf{R}^n)$ le champ mesurable de matrices positives de $g_i$ par rapport à la structure standard, de telle sorte que la norme sur $\Omega^k(U_i, g_i)$ est donné par:

$$\left\{ \int_{U_i} \langle \lambda^k(A^{-1})\omega, \omega \rangle \det(A)^{\frac{1}{2}} dx \right\}^{\frac{1}{2}}.$$

Si $a \in M_n(\mathbf{R})$ est une matrice positive et $\geq 1$, on a toujours

$$\lambda^k(a^{-2}) \det(a) \leq \|a\|^{n-k} \quad \text{et} \quad \lambda^k(a^2) \det(a)^{-1} \leq \|a\|^k.$$

Par conséquent, $\lambda^k(A^{-1}) \det(A)^{\frac{1}{2}} \in L^{\frac{p}{n-k}}$, et $\lambda^k(A) \det(A)^{-\frac{1}{2}} \in L^{p/k}$, ce qui montre, encore par l'inégalité de Hölder, que nous avons les inclusions:

$$L^{p_k} \subset \Omega^k(V, g) \subset L^{q_k}$$

avec $p_k = \frac{2p}{p+k-n}$ et $q_k = \frac{2p}{p+k}$. Si $n = 2m+1$, les relations de la définition 1.1 sont donc satisfaites dès que $\frac{p-m-1}{2p} + \frac{1}{2n} > \frac{p+m+1}{2p}$, c'est à dire dès que:

$$p > n(m+1) = \frac{n(n+1)}{2}.$$

Si $n = 2m$, on doit avoir $\frac{p-m}{2p} + \frac{1}{n} > \frac{p+m+1}{2p}$, ce qui donne:

$$p > \frac{n(2m+1)}{2} = \frac{n(n+1)}{2}.$$

## 4. Conséquences en K-théorie

Si $V$ est une variété de dimension vérifiant les conditions du théorème 2.7, le cycle $(\Omega^*(V,g), D)$ associé à $g \in \mathcal{R}(\mathcal{O})$ définit alors une classe $\Sigma_V \in K_n(V)$ ($n = \dim V$) indépendante des choix de $\mathcal{O}$ et de $g$. Si $V$ est de dimension paire, le cycle de Fredholm associé est $(\Omega^m(V,g), F)$ avec $F = (dd^* - d^*d)(1 + dd^* + d^*d)^{-1}$. En effet, les opérateurs $F^2 - 1$ et $F - F^*$ sont compacts et $F\tau + \tau F = 0$; pour tout $\varphi \in C(V)$, on a $[F, \varphi] \in \mathcal{K}$. Dans le cas impair, le cycle de Fredholm associé est donné par $(\Omega^m(V,g), F)$ avec $F = (\tau d + 1)(1 + (\tau d)^2)^{-\frac{1}{2}}$.

Cet élément $\Sigma_V$ ainsi construit représente la classe de l'opérateur de signature usuel d'une variété riemannienne $C^\infty$. Cela est démontré dans ([5, Th. 3.2]) si $\dim V$ est paire. Notons que si $V$ est une variété quasi-conforme de dimension paire, le module de Fredholm borné $(\Omega^m(V,g), F)$ est égal, modulo les compacts, au module construit dans [5].

Dans le cas impair, en supposant par exemple que $m$ est pair, on a une décomposition orthogonale de l'identité de $\Omega^{(0)} = \sum \Omega^{2k}$, $1 = e + f + g$, avec $\operatorname{im} e = (\operatorname{supp} \tau d) \cap \Omega^m(V,g)$, $\operatorname{im} f = \sum_{2k>m} \Omega^{2k}(V,g)$ et dans laquelle l'opérateur de signature $D = \tau d - d\tau$ a la forme suivante (cf. [2, Prop. 4.20]):

$$D = \begin{pmatrix} \tau d & 0 & 0 \\ 0 & 0 & A^* \\ 0 & A & 0 \end{pmatrix}.$$

L'élément de K-homologie associé est alors représenté par $e_1 + f_1$, où $e_1$ est le projecteur spectral positif de la restriction de $\tau d$ à $\operatorname{im} e$ et $f_1$ celui de la restriction à $\operatorname{im} f \oplus \operatorname{im} g$. On a alors $\operatorname{im} f_1 = \text{graphe } U$ où $U$ est la phase de $A$ et qui vérifie $[\varphi, U] \in \mathcal{K}$ pour $\varphi \in C(V)$. La famille à un paramètre de projections donnée par $e_1 + f_t$, où $f_t$ est le projecteur orthogonal de $\Omega^{(0)}$ sur graphe$(tU)$, est une homotopie opératorielle entre $e_1 + f_1$ et $e_1 + f_2$, où $f_2$ est un projecteur tel que $f_2 \subset f$ et $\dim(f - f_2) < +\infty$. Le projecteur $e_1 + f_2$ est alors égal modulo un projecteur de rang fini à $e_1 \oplus f_3$, où $f_3$ est une projection hermitienne d'un $C(V)$-module orthogonal á $\Omega^m(V,g)$.

Soit ch : $K_i(V) \otimes \mathbf{Q} \to H_*(V, \mathbf{Q})$ le caractère de Chern et $\mathcal{P} : H^*(V, \mathbf{Q}) \to H_*(V, \mathbf{Q})$ la dualité de Poincaré et $\mathcal{L}(V)$ le polynôme de Hirzebruch en les classes de Pontryagin rationnelles de $V$. Les résultats présents permettent de construire comme dans [8] une $K$-orientation pour le microfibré tangent à de telles variétés et nous obtenons alors:

THÉORÈME 4.1.    *Soit $V$ une variété compacte orientée de dimension $n$, ou bien quasi-conforme, ou bien dérivable à l'ordre $p > \frac{n(n+1)}{2}$. Pour tout recouvrement fini par des ouverts de carte $\mathcal{O}$ et pour tout $g \in \mathcal{R}(\mathcal{O})$, le cycle non-borné $(\Omega(V, g), D)$ détermine un élément $\Sigma_V \in K_n(V)$ qui ne dépend que de la structure ou bien quasi-conforme, ou bien dérivable à l'ordre $p$ de $V$. La classe $\Sigma_V$ est une base de $K_i(V) \otimes \mathbf{Z}[\frac{1}{2}]$ en tant que module sur $K^i(V)$ et on a:*

$$\mathrm{ch}\, \Sigma_V = \mathcal{P}(\mathcal{L}(V)).$$

UNIVERSITÉ P. ET M. CURIE (B.P. 191), 75 252 PARIS CEDEX 05, FRANCE
*E-mail address*: hilsum@math.jussieu.fr